\documentstyle{article}
\begin{document}

\newtheorem{lemma}{Lemma}
\newtheorem{theorem}{Theorem}
\newtheorem{definition}{Definition}
\newtheorem{corollary}{Corollary}
\newtheorem{proposition}{Proposition}

\def\R{{\bf R}}
\def\C{{\bf C}}
\def\Z{{\bf Z}}
\def\bc{\vskip3mm\begin{center}}
\def\ec{\end{center}\vskip3mm}
\def\be{\begin{equation}}
\def\ee{\end{equation}}
\def\bt{\begin{theorem}}
\def\et{\end{theorem}}
\def\bl{\begin{lemma}}
\def\el{\end{lemma}}
\def\bd{\begin{definition}}
\def\ed{\end{definition}}
\def\bcor{\begin{corollary}}
\def\ecor{\end{corollary}}
\def\bp{\begin{proposition}}
\def\ep{\end{proposition}}

\bc
{\large The Weierstrass representation of closed surfaces in $\R^3$}
\ec
\bc
{\large Iskander A. Taimanov
\footnote{Institute of Mathematics,
630090 Novosibirsk, Russia;
e-mail : taimanov@math.nsc.ru}}
\ec

\bc
{\bf \S 1. Introduction}
\ec

The present article 
is a sequel to \cite{T1,T2}.
The results presented here extend onto general surfaces 
the results obtained in \cite{T2} for surfaces of revolution
and were exposed in a lot of talks of 
the author during the last year being  at the end announced in \cite{T3}.
These are

-- a construction of a global Weierstrass representation for an
arbitrary closed oriented surface, of genus $g \geq 1$, immersed into
$\R^3$ (\S 2);

-- a construction of a Weierstrass spectrum for a torus immersed into
$\R^3$ and discussing its geometric properties (\S 3);

-- a construction of finite-zone surfaces and 
finite-zone solutions to the modified Novikov--Veselov equations 
(\S\S 4-5).

In \S 6 we discuss a relation of these constructions
to the Willmore conjecture.

\bc
{\bf \S 2. The Weierstrass representation}
\ec

{\bf 2.1. The local Weierstrass representation.}

The local Weierstrass representation of surfaces immersed into $\R^3$ is based
on the following two lemmas.

\bl
Let $W$ be a simply connected domain in $\C$ and
let a vector function $\psi = (\psi_1, \psi_2): W \rightarrow \C^2$ satisfy
the equation
\be
{\cal D}\psi = 0
\label{1}
\ee
where
\be
{\cal D} =
\left(
\begin{array}{cc}
0 & \partial \\
- \bar{\partial} & 0
\end{array}
\right)
+
\left(
\begin{array}{cc}
U & 0 \\
0 & U
\end{array}
\right)
\label{2}
\ee
and a function $U(z,\bar{z})$ is real-valued. Then the following formulas
$$
X^1(z,\bar{z}) = \frac{i}{2} \int_{\gamma}
\left( (\bar{\psi}_2^2 + \psi_1^2) dz' -
(\bar{\psi}_1^2 + \psi_2^2) d\bar{z}' \right),
$$
\be
X^2(z,\bar{z}) = \frac{1}{2}\int_{\gamma}
\left(
(\bar{\psi}_2^2 - \psi_1^2) dz' -
(\bar{\psi}_1^2 - \psi_2^2) d\bar{z}'
\right),
\label{3}
\ee
$$
X^3(z,\bar{z}) = \int_{\gamma}
(\psi_1 \bar{\psi}_2 dz' + \bar{\psi}_1 \psi_2 d\bar{z}')
$$
define an immersion of $W$  into $\R^3$. Moreover the induced metric
takes the form
$$
D(z,\bar{z})^2 dz d\bar{z},
\ \mbox{where} \  
D(z,\bar{z}) = (|\psi_1(z,\bar{z})|^2 + |\psi_2(z,\bar{z})|^2),
$$
the Gauss curvature is
$$
K(z,\bar{z}) = -\frac{4}{D(z,\bar{z})^2}
\partial \bar{\partial}
\log D(z,\bar{z}),
$$
and the mean curvature equals
\be
H(z,\bar{z}) = 2 \frac{U(z,\bar{z})}{D(z,\bar{z})}.
\label{4}
\ee
\el

The immersion $W \rightarrow \R^3$ is constructed as follows.
Take $z_0 \in W$ and map it into $0 \in \R^3$.
After that define $X^i(z,\bar{z})$ by the value of the integral (\ref{3}) 
taken over an arbitrary path, in $W$, connecting $z_0$ and $z$.

For $U = 0$, this is the classical Weierstrass representation 
of minimal surfaces. 

After this lemma is formulated it is not difficult to prove it.
The condition (\ref{1}) together with the reality of $U$ imply that the
integrands in (\ref{3}) are closed forms. Hence (\ref{3}) defines an immersion
$W$ into $\R^3$ for which it is easy to compute all characteristics.

It seems that this lemma belongs to Eisenhart (\cite{Eis})
who had written (\ref{1}) as a condition only for $\psi_1$:
$$
L\psi_1 = 0,
$$
in terms of the second order scalar differential operator
\be
L = \partial \bar{\partial} -
\frac{\partial U}{U}\bar{\partial}  + U^2.
\label{5}
\ee

In late 80's Abresch had derived the Weierstrass representation for
constant mean curvature surfaces using the operator (\ref{2})
and considered global representations for constructing
explicit formulae for such surfaces of genus $g \geq 2$. Since this program
is not realized until now, these results were only presented on
the Luminy conference (1989).

Bobenko had shown that considering $\R^3$ 
as the space of imaginary quaternions
some identities for the moving frame are written in terms of the Dirac
equation for quaternion-valued functions,
but constructing of general surfaces in
terms of eigenfunctions of the Dirac operator and the
globalization of this construction were not discussed by him (\cite{Bo}).  

In \cite{Kon} Konopelchenko considered a ``square root'' of 
(\ref{5}), i.e., the first order matrix differential operator (\ref{2}),
for definition of deformations of surfaces given by the formulas
(\ref{3}) via the modified Novikov--Veselov equations. 

This replacement of the two-dimensional Schr\"odinger operator with 
a potential and a magnetic field (\ref{5}) by the Dirac operator (\ref{2}) 
with the potential $U$ makes a representation more 
comfortable for applications.

\bl
Let $W$ be a domain in $\C$ and let $X: W \rightarrow \R^3$ be a
conformal immersion of $W$ into $\R^3$:
$z \rightarrow X(z,\bar{z}) = (X^1(z,\bar{z}), X^2(z,\bar{z}),
X^3(z,\bar{z}))$. Assume that 
\be
\frac{\partial X^3}{\partial z} \neq 0
\label{6}
\ee
near $z_0 \in W$. Then near $z_0$ the functions
\be
\psi_1(z,\bar{z}) =
\sqrt{-\partial \Phi(z,\bar{z})},\ \ \
\psi_2(z,\bar{z}) =
\sqrt{\bar{\partial} \Phi(z,\bar{z})},
\label{7}
\ee
with
$$
\Phi(z,\bar{z}) = X^2(z,\bar{z}) + i X^1(z,\bar{z}),
$$
satisfy (\ref{1})
with $U(z,\bar{z}) = H(z,\bar{z}) D(z,\bar{z})/2$,
where $H$ is the mean curvature and $D^2 dz d\bar{z}$ is the metric 
of the surface $X(W) \subset \R^3$.
\el

We explain the proof in brief following \cite{T1}.

A condition that an immersion is conformal is written as
\be
\left(\frac{\partial X^1}{\partial z}\right)^2 +
\left(\frac{\partial X^2}{\partial z}\right)^2 +
\left(\frac{\partial X^3}{\partial z}\right)^2 = 0.
\label{8}
\ee
The quadric $Q = \{x_1^2 + x_2^2 + x_3^2 = 0\} \subset \C P^2$
is diffeomorphic to the Grassmann manifold of two-dimensional oriented 
linear subspaces of $\R^3$ and for conformally immersed surface $X(W)$ 
the mapping $z \in W \rightarrow (X^1_z : X^2_z : X^3_z) \in Q$
is the Gauss map.

The Gauss map for the surface 
(\ref{3}) is related to $\psi$ by (\ref{7}).
Consider a conformal immersion satisfying (\ref{6}) and by (\ref{7}) 
construct near $z_0 \in W$ a vector function $\psi$.
It follows from (\ref{6}) and (\ref{8}) that near $z_0$
the radicands in (\ref{7}) do not vanish.

Take the following bases for the tangent planes to $X(W)$:
$$
e_1 = \frac{1}{D} \frac{\partial X}{\partial x}, \ \ \
e_2 = \frac{1}{D} \frac{\partial X}{\partial y},
$$
where $D^2 dz d\bar{z}$ is the induced metric.  Put $e_3 = e_1 \times e_2$.
The derivation formulas take the form
$$
\frac{\partial^2 X}{\partial x^2} =
\frac{\partial D}{\partial x}e_1 -
\frac{\partial D}{\partial y}e_2 + D^2h_{11}e_3,
$$
\be
\frac{\partial^2 X}{\partial x \,\partial y} =
\frac{\partial D}{\partial y}e_1 + \frac{\partial D}{\partial x}e_2 +
D^2h_{12}e_3,
\label{9}
\ee
$$
\frac{\partial^2 X}{\partial y^2} =
-\frac{\partial D}{\partial x}e_1 + \frac{\partial D}{\partial y}e_2 +
D^2h_{22}e_3,
$$
where  $h_{ij}$ is the second fundamental form.
 
Now Lemma 2 is proved by substitution of (\ref{9}) into the expressions for
$\bar{\partial} \psi_1$ and $\partial \psi_2$ and by consequent
straightforward computations.

\bd
A representation of a surface $\Sigma$, immersed into $\R^3$,
by the formulas (\ref{3}) is called a Weierstrass representation.

The function $U(z,\bar{z})$ of the form (\ref{4})
is called the potential of the surface $\Sigma$ with the distinguished
conformal parameter $z$, or the potential of the Weierstrass representation 
of $\Sigma$.
\ed

Lemma 2 immediately implies (\cite{T1})

\bl
Every regular $C^2$-surface immersed into $\R^3$ near every its point 
admits a Weierstrass representation.
\el

{\bf 2.2. The global Weierstrass representation.}

For defining a Weierstrass representation globally for the whole surface 
$\Sigma$ it needs to correctly define a bundle generated by $\psi$ over 
$\Sigma$ and an operator ${\cal D}$ acting in this bundle. This problem was
solved in \cite{T1} and we recall this solution.

We consider closed oriented surfaces of genus $g \geq 1$.

By the Riemann uniformization theorem, every torus is conformally equivalent 
a flat torus $\C/\Lambda$ with $\Lambda$ a lattice of rank $2$.

\bd
A torus $\Sigma$, immersed into $\R^3$, possesses a global Weierstrass 
representation if there exist a real potential $U$ and  functions
$\psi_1$ and $\psi_2$ defined on the universal covering of $\Sigma$, i.e.,
on $\C$, such that 

1)
\be
\left\{
\begin{array}{l}
U(z + \gamma) = U(z), \\
\psi_j(z+\gamma) = \varepsilon(\gamma) \psi_j(z), \\
\varepsilon(\gamma)=\pm 1
\end{array}
\right.
\label{10}
\ee
for $z \in \C$ and $\gamma \in \Lambda$;

2) the vector function $\psi$ satisfies (\ref{1}) and  
for a suitable choice of coordinates in $\R^3$ defines by (\ref{3})
an immersion of $\Sigma$.
\ed

For a sphere $\Sigma$ with $g >1$ handles the uniformization theorem 
tells that $\Sigma$ is conformally equivalent to a quotient 
${\cal H}/\Lambda$ with ${\cal H}$ the Lobachevskii upper-half plane and
$\Lambda$ a lattice in $PSL(2,\R)$, the isometry group of ${\cal H}$.

\bd
A sphere $\Sigma$ with $g \,(>1)$ handles, immersed into $\R^3$,
possesses a global Weierstrass representation if there exist a real
potential $U$ and functions $\psi_1$ and $\psi_2$, defined on the universal
covering of $\Sigma$, i.e., on ${\cal H}$, such that

1)
\be
\left\{
\begin{array}{l}
U(\gamma(z)) = |cz +d|^2 U(z), \\
\psi_1(\gamma(z)) = (cz+d) \psi_1(z), \\
\psi_2(\gamma(z)) = (c\bar{z}+d) \psi_2(z)
\end{array}
\right.
\label{11}
\ee
for $z \in {\cal H}$ and $\gamma \in \Lambda$, represented by the matrix
$$
\left(
\begin{array}{cc}
a & b \\
c & d
\end{array}
\right),
\ \ \ \
a,b,c,d \in \R, \ \
ad-bc=1;
$$

2) the vector function $\psi$ satisfies (\ref{1}) and 
for a suitable choice of coordinates in $\R^3$ 
defines by (\ref{3}) an immersion of $\Sigma$.
\ed

These definitions are quite natural and based on the following statement.

\bl
If a local Weierstrass representation of a closed oriented surface $\Sigma$
of genus $g \geq 1$ is smoothly continued onto the whole universal
covering $\tilde{\Sigma}$ then $U$ and 
$\psi$ satisfy (\ref{10}) for $g=1$ and (\ref{11}) for $g\geq 2$.  
\el

This lemma follows from the transformation rules for the mean curvature and 
the metric under changing of coordinates. Now we arrive 
at the following result(\cite{T1}).

\bt
If a closed oriented surface $\Sigma$ of genus $\geq 1$ possesses a global
Weierstrass representation then $\psi$ is a global section
of a spinor bundle over the constant curvature surface $\Sigma_0$ 
conformally equivalent to $\Sigma$ and the Dirac operator (\ref{2}) 
acts on this bundle.
\et

The following theorem demonstrates an importance of this representation.

\bt
\footnote{This theorem also holds for spheres but we will consider
this case separately.}
Every $C^3$-regular closed oriented surface $\Sigma$ of genus
$g \geq 1$, immersed into $\R^3$, possesses 
a global Weierstrass representation.
\et

A proof of Theorem 2.

By Lemma 3, the functions (\ref{7}) define a local Weierstrass representation 
near every point with $\partial \Phi \neq 0$.
For proving an existence of a global representation it suffices to correctly 
smoothly extend these functions onto neighborhoods of points with 
$\partial \Phi = 0$. By (\ref{8}), the last condition is 
equivalent to $X^3_z = 0$.  

For any compact immersed surface 
$\Sigma \subset \R^3$ there exist coordinates in $\R^3$ such 
that all the critical points of the ``height function'', 
a distance to the  plane
$X^3 =0$, are nondegenerate in Morse's sense, i.e., 
at each of them the matrix of second derivatives of 
$X^3(z,\bar{z})$ is nondegenerate (\cite{Milnor}).
Take such coordinate system for the surface in study $\Sigma$.

Outside critical points define $\psi$ by (\ref{7}) and show that these 
functions do not ramify at critical points.

Let $q$ be a critical point of $X^3$ and let $z$ be 
a conformal parameter near $q$ such that $z(q)=0$.
Since this point is nondegenerate, we have
\be
X^3(z,\bar{z}) = \alpha z^2 +
\bar{\alpha}\bar{z}^2 + \beta z\bar{z} + O(|z|^3)
\label{12}
\ee
with $|\alpha|+|\beta| \neq 0$.
Show that (\ref{7}) correctly defines $\psi_1$ near $q$ and notice
that these conversations work also for $\psi_2$. 

Since $\partial \Phi (q) = 0$, we have
\be
\frac{\partial X^1(q)}{\partial z} =
i\frac{\partial X^2(q)}{\partial z}.
\label{13}
\ee

Differentiating the left-hand side of (\ref{8}) by $z$, we infer from
(\ref{12}) and (\ref{13}) that
\be
(\partial \Phi)_z (q) = 
\frac{\partial^2 X^2(q)}{\partial z^2} +
i\frac{\partial^2 X^1(q)}{\partial z^2} = 0.
\label{14}
\ee
Analogously differentiating it by $\bar{z}$, we derive
$$
(\partial \Phi)_{\bar{z}} (q) =
\frac{\partial^2 X^2(q)}{\partial z \, \partial \bar{z}} +
i\frac{\partial^2 X^1(q)}{\partial z\, \partial \bar{z}}
= 0.
$$

Assume that 
$\alpha \neq 0$.
Differentiating the left-hand side of (\ref{8})
twice by $z$ and taking  (\ref{12}), (\ref{13}) and (\ref{14})
into account, we obtain
$$
\left(\frac{\partial^2 X^3(q)}{\partial z^2}\right)^2 +
\frac{\partial X^2(q)}{\partial z}
\left(\frac{\partial^3 X^2(q)}{\partial z^3} + 
i\frac{\partial^3 X^1(q)}{\partial z^3}\right) = 0
$$
and, since the surface is regular, $X^2_z(q) \neq 0$. 
Hence 
\be
(\partial \Phi)_{zz}(q) \neq 0.
\label{15}
\ee

If $\alpha =0$, then 
$\beta \neq 0$ and, differentiating the left-hand side of
(\ref{8}) twice by $\bar{z}$,
we obtain
\be
(\partial \Phi)_{\bar{z}\bar{z}}(q) \neq 0.
\label{16}
\ee

From (\ref{15}) and (\ref{16}) we conclude that
$\partial \Phi$ has a double zero at $q$ and if, for instance, 
the inequality (\ref{15}) holds then
$\partial \Phi = z^2 \cdot f(z,\bar{z})$ where the branches of
$\sqrt{f}$ do not ramify at $q$.
Hence the branches of $\sqrt{-\partial \Phi}$ do not ramify at $q$.

Since we consider an arbitrary critical point of the function $X^3$,
each branch of $\sqrt{-\partial \Phi}$
is correctly defined and has no ramifications.
Therefore on the universal covering of $\Sigma$ the vector function $\psi$ is
correctly defined by (\ref{7}).

This proves Theorem 2.

We arrive at the following conclusion

{\sl Zero-eigenfunctions of Dirac operators (\ref{2}), with 
real potentials, on spinor bundles (\ref{10}) and (\ref{11})  
over closed oriented constant curvature surfaces $\Sigma$ of genus
$\geq 1$ are in one-to-one correspondence with immersions of their universal 
coverings into $\R^3$ with Gauss maps descending through $\Sigma$.}

\bc
{\bf \S 3. The spectrum of Weierstrass representation}
\ec

Let $\Sigma$ be a torus immersed into $\R^3$ and conformally equivalent
to a flat torus $\C/\Lambda$.
Assume that it possesses a global Weierstrass representation
which is, by Theorem 2, valid if it is $C^3$-regular :

1) There exists a function $U(z): \C \rightarrow \R$,
the potential of $\Sigma$,
such that
$U(z + \gamma) = U(z)$ for $z \in \C, \gamma \in \Lambda$;

2) There exist $\varepsilon_1$ and $\varepsilon_2$ such that
$\varepsilon_j \in \{0,1\}$ and they determine a spinor bundle over
$\C / \Lambda$;

3) There exists a function $\psi$ such that ${\cal D}\psi = 0$
and it satisfies the periodicity conditions (\ref{10}).

By Lemma 1, every solution to (\ref{1}) defines an immersion of a surface 
into $\R^3$. We consider a ``linear basis'' for the family of 
such surfaces, i.e., the set of Floquet functions of the 
operator (\ref{2}).

\bd
A function $\psi: \C \rightarrow \C$ is called a Floquet function
(with zero eigenvalue) of the operator ${\cal D}$ (\ref{2})
with the quasimomenta $(k_1, k_2)$ if ${\cal D} \psi = 0$ and
\be
\psi(z + \gamma) =
\exp{(2\pi i
(\mbox{Re}\, \gamma \cdot k_1 + \mbox{Im}\, \gamma \cdot k_2)
)}\psi(z)
\label{17}
\ee
for $\gamma \in \Lambda$.
\ed

Notice that
a function $\psi$ satisfying (\ref{17}) has the form
\be
\psi(z) = \exp{(2\pi i(xk_1 + yk_2))}\varphi(z)
\label{18}
\ee
where $\varphi(z)$ is periodic with respect to $\Lambda$.

The following lemma is clear.

\bl
A function $\psi$ of the form (\ref{18}) satisfies the equation
\be
\left[
\left(
\begin{array}{cc}
0 & \partial \\
- \bar{\partial} & 0
\end{array}
\right)
+
\left(
\begin{array}{cc}
U & 0 \\
0 & U
\end{array}
\right)
\right] \psi = \lambda\psi
\label{19}
\ee
if and only if
\be
{\cal D}_{k}\varphi = \lambda\varphi
\label{20}
\ee
where
$$
{\cal D}_{k} =
\left(
\begin{array}{cc}
0 & \partial \\
- \bar{\partial} & 0
\end{array}
\right)
+
\left(
\begin{array}{cc}
U &  \pi(k_2 + i k_1)
\\
\pi(k_2 - i k_1) & U
\end{array}
\right).
$$
\el

Take a constant $C$ such that the operator 
$$
{\cal A} = 
\left[
\left(
\begin{array}{cc}
0 & \partial \\
- \bar{\partial} & 0
\end{array}
\right)
+
\left(
\begin{array}{cc}
C & 0 \\
0 & C
\end{array}
\right)
\right]
$$
is invertible on $L_2(\C/\Lambda)$.
An  existence of such constant is easily verified by 
using the Fourier transform.
Then the equation (\ref{20}) possesses a solution from $L_2(\C/\Lambda)$
if and only if the equation
\be
\left[ 1 +
\left(
\begin{array}{cc}
U - (C + \lambda) & \pi(k_2 +  i k_1) \\
\pi(k_2 - i k_1) & U - (C + \lambda)
\end{array}
\right)
\left(
\begin{array}{cc}
C & \partial \\
- \bar{\partial} & C
\end{array}
\right)^{-1}
\right] \xi = 0
\label{21}
\ee
is solvable in $L_2(\C/\Lambda)$.

Since ${\cal A}^{-1}$ increases a smoothness, by the
Sobolev embedding theorem, ${\cal A}^{-1}$ is compact.
Since $U$ is continuous, the operator of 
multiplication by
$$
\left( \begin{array}{cc} U - (C + \lambda) & \pi(k_2 +
i k_1) \\ \pi(k_2 - i k_1) & U - (C + \lambda)
\end{array} \right).
$$
is bounded. Hence, we have

\bl
A pencil of operators
$$
{\cal D}_{k} \circ
\left(
\begin{array}{cc}
C & \partial \\
- \bar{\partial} & C
\end{array}
\right)^{-1} -
1  =
$$
\be
\left(
\begin{array}{cc}
U - (C + \lambda) & \pi(k_2 + i k_1) \\
\pi(k_2 - i k_1) & U - (C + \lambda)
\end{array}
\right)
\left(
\begin{array}{cc}
C & \partial \\
- \bar{\partial} & C
\end{array}
\right)^{-1} 
\label{22}
\ee
is polynomial in $k_1, k_2$ and $\lambda$
and consists of compact operators from 
$L_2(\C/\Lambda)$ to $L_2(\C/\Lambda)$.
\el

Now by using the polynomial Fredholm alternative 
(\cite{Kel}) we obtain

\bt 
There exist analytic subsets $\hat{\Gamma}_U \subset \C^3$
and $\Gamma_U \subset \C^2$ of positive codimensions such that

1) the equation (\ref{21}) is solvable in $L_2(\C/\Lambda)$ if
and only if $(k_1, k_2, \lambda) \in \hat{\Gamma}_U$ ;

2) $\Gamma_U = \hat{\Gamma}_U \cap \{\lambda=0\}$.
\et

These subsets are called the Floquet spectrum and the zero 
Floquet spectrum of ${\cal D}$, respectively.

A proof of Theorem 3.

By the Keldysh theorem (\cite{Kel}), if $A_{\mu}$ is a polynomial, in
$\mu \in \C^n$,  pencil of compact operators, then the set of
$\mu$, for which the equation $(1 + A_{\mu})\xi = 0$ is solvable,
forms an analytic subvariety in $\C^n$.

We are left to prove that these subsets have positive codimensions.

Let $k_1 = \lambda = 0$.
Then (\ref{21}) is equivalent to
\be
\left(
\left(
\begin{array}{cc}
C & \partial \\
-\bar{\partial} & C
\end{array}
\right)
+
\left(
\begin{array}{cc}
0 & \pi k_2 \\
\pi k_2 & 0
\end{array}
\right)
\left(
\begin{array}{cc}
1 & \frac{U-C}{\pi k_2} \\
\frac{U-C}{\pi k_2} & 1
\end{array}
\right)
\right) \xi = 0
\label{23}
\ee
and as $|k_2| \rightarrow \infty$
it degenerates into 
$$
\left(
\left(
\begin{array}{cc}
0 & \partial \\
-\bar{\partial} & 0
\end{array}
\right)
+
\left(
\begin{array}{cc}
0 & \pi k_2 \\
\pi k_2 & 0
\end{array}
\right)
\right) \xi = 0.
$$
From the last equation by methods of perturbation theory it is inferred that
there exists $k_2$ (with sufficiently large $|k_2|$) for which the equation
(\ref{23}) is not solvable in $L_2(\C/\Lambda)$. This proves Theorem 3.

Kuchment had strengthened the Keldysh theorem for some operators:

{\sl There exists an entire function $Z: \C^3 \rightarrow \C$
such that the equation (\ref{21}) is solvable in $L_2(\C/\Lambda)$
if and only if $Z(k_1,k_2,\lambda)=0$.}

Roughly speaking, $Z$
is a regularized determinant of the pencil
(\ref{22}).
In \cite{Ku} this is proved for scalar even order elliptic operators 
but his reasonings work also for the operator (\ref{2}).

It follows from (\ref{18}) that $\hat{\Gamma}_U$ and $\Gamma_U$ are 
invariant under the action of the dual lattice $\Lambda^*$:
\be
k_1 \rightarrow k_1 + \mbox{Re}\,\gamma^*,   
k_2 \rightarrow k_2 + \mbox{Im}\,\gamma^*,  
\gamma^*\in \Lambda^*.
\label{24}
\ee
Remind that $\Lambda^*$ consists of $\gamma^* \in \C$ such that
$(\gamma,\gamma^*) = \mbox{Re}\,\gamma \cdot
\mbox{Re}\,\gamma^* + \mbox{Im}\,\gamma \cdot \mbox{Im}\,\gamma^* 
\in \Z$ for all $\gamma \in \Lambda$.
It is also clear that the definitions of
$\hat{\Gamma}_U$ and $\Gamma_U$ are independent on a choice of
a conformal parameter $z$ on the torus.

\bd
A complex surface $W_{\Sigma}$ defined as the quotient space of $\Gamma_U$ 
for the action (\ref{24}) is called the Weierstrass spectrum of the torus 
$\Sigma$.

The genus of the normalization of $W_{\Sigma}$ is called the spectral genus
of $\Sigma$.
\ed

This spectrum was introduced by the author who had pointed out its relation 
to the Willmore functional, i.e., an integral of the 
squared mean curvature, as the simplest Kruskal integral. 
For tori of revolution this was analyzed in details in \cite{T2}. 
  
More effective definition of the Floquet spectrum for multidimensional
operators is given by perturbation theory. For the Schr\"odinger
operator and for $\partial_y - \partial^2_x$ 
this was done by Krichever (see \cite{Kr}, where 
an approach to define spectral curves by using
the Keldysh theorem is also mentioned with a reference to our unpublished 
paper).  

The Willmore functional ${\cal W}$ is conformally invariant in the following
sense. Let $g: \R^3 \rightarrow \R^3$ be a conformal transformation
saving a closed surface $\Sigma$ in a compact domain, then 
${\cal W}(\Sigma) = {\cal W}(g(\Sigma))$.
The author had conjectured that the whole Weierstrass spectrum is 
conformally invariant in this sense.
Quite soon after its formulation two different proofs of this conjecture 
had been  obtained.

Pinkall using methods of (\cite{KPP}) 
had written the spectral problem ${\cal D}\psi = 0$
in conformally invariant terms which implies the conformal invariance of 
$W_\Sigma$.

Grinevich and Schmidt had shown that,
since the transformation formulas, for the potential, corresponding to
infinitesimal conformal transformations of a surface are quadratic in
$\psi$ (it was noticed for tori in revolution in \cite{T2}), 
by an analog of the Melnikov theorem for periodic operators, 
$\Gamma_U$ is conformally invariant (\cite{GS}).

A relation of the Weierstrass spectrum to the spectral curves of special 
soliton tori (constant mean curvature, Willmore, see \cite{FPPS,H,PS}) 
will be discussed elsewhere.

\bc
{\bf \S 4. Finite-zone planes and tori}
\ec

In the next two paragraphs we somewhere only sketch proofs which are
usual for the finite-zone theory (\cite{DKN,Kr0}).
Moreover for the one-dimensional limit of (\ref{2}) a derivation of 
the theta formulas is exposed in \cite{Prev} and  symmetries of $\Gamma$
are discussed in \cite{D}.

Consider more general operator
\be
L=
\left(
\begin{array}{cc}
0 & \partial \\
-\bar{\partial} & 0
\end{array}
\right)
+
\left(
\begin{array}{cc}
U & 0 \\
0 & V
\end{array}
\right)
\label{25}
\ee
with $U$ and $V$ periodic with respect to a rank two lattice
$\Lambda \subset \C$. 
For this operator an analog of Theorem 3 holds
and $L$ is called finite-zone 
(on the zero energy level) if the normalization of its zero Floquet spectrum 
is a compact Riemann surface with two points removed.
But we will call $L$ finite-zone 
if it is as follows.
\footnote{It is clear that the equivalence of these definitions for 
smooth potentials can be justified by perturbation theory
(\cite{Kr0}).}

\bp
Let $\Gamma$ be a compact Riemann surface of genus $g$,
$\infty_{\pm}$ be a pair of distinct points on $\Gamma$,
$k^{-1}_{\pm}$ be local parameters near these points such that
$k^{-1}_{\pm}(\infty_{\pm}) =0$, and $D$ be a nonspecial effective
divisor of degree $g+1$ on $\Gamma \setminus \{\infty_{\pm}\}$,
i.e., $D=P_1 + \dots + P_{g+1}$ with $P_i \in \Gamma 
\setminus \{\infty_{\pm}\}$. Then 

1. There exists a unique vector function
$\psi(z,\bar{z},P) = (\psi_1, \psi_2)$, with $z \in \C$, such that  
$\psi$ is meromorphic in $P$ on $\Gamma \setminus \{\infty_{\pm}\}$
and has poles only in $D$,

\be
\psi = \exp{(k_+z)} 
\left[
\left( \begin{array}{c} 1 \\ 0  \end{array} \right)
+
\left( \begin{array}{c} \xi^+_{11}/k_+ \\ \xi^+_{21}/k_+ \end{array} 
\right)
+ O(k_+^{-2}) 
\right] \ \ \ \ 
\mbox{as $P \rightarrow \infty_+$},
\label{26}
\ee
and
\be
\psi = \exp{(k_- \bar{z})} 
\left[
\left( \begin{array}{c} 0 \\ 1  \end{array} \right)
+
\left( \begin{array}{c} \xi^-_{11}/k_- \\ \xi^-_{21}/k_- \end{array} 
\right)
+ O(k_-^{-2}) 
\right] \ \ \ \ 
\mbox{as $P \rightarrow \infty_-$}.
\label{27}
\ee

2. Moreover there exists a unique operator $L$ of the form (\ref{25})
such that $L \psi = 0$.
The potentials of $L$ are as follows
\be
U = -\xi^+_{21}, \ \ \
V = \xi^-_{11}.
\label{28}
\ee
\ep  

This proposition is a particular case of the general theorem on uniqueness
of the Baker-Akhieser function (\cite{Kr0}).  

Fix a basis $\alpha_1,\dots,\alpha_g,\beta_1,\dots,\beta_g$ for $H_1(\Gamma)$
such that the intersection form is 
$$
\alpha_j \circ \beta_k = \delta_{jk}, \ \ 
\alpha_j \circ \alpha_k = \beta_j \circ \beta_k = 0.
$$
To this basis corresponds a unique basis of holomorphic differentials
$\omega_1, \dots, \omega_g$ normalized by the condition
$$
\int_{\alpha_k} \omega_j = \delta_{jk}.
$$
Define now the period matrix $\Omega$ and the theta function
of $\Gamma$ as follows
$$
\Omega_{jk} = \int_{\beta_k} \omega_j,
$$
$$
\vartheta(u) = \sum_{N \in \Z^g} 
\exp{\pi i ((\Omega N,N) + 2(N,u))},
$$
where $u \in \C^g$.
Fixing a point $P_0 \in \Gamma \setminus \{\infty_{\pm}\}$,
define also the Abel map from $\Gamma$ into its Jacobian variety
$J(\Gamma) = \C^g / \{M + \Omega N : M,N \in \Z^g\}$:
$$
A(P) = (\int_{P_0}^P \alpha_1, \dots, \int_{P_0}^P \alpha_g).
$$

Denote by $\eta_l^{\pm}$ a unique meromorphic differential having a single
pole at $\infty_{\pm}$ with the Laurent part $d k_{\pm}^l$ and 
normalized by the condition
$$
\int_{\alpha_j} \eta_l^{\pm} = 0 \ \ \ \
\mbox{for $j =1,\dots,g$.}
$$
To every $\eta^{\pm}_l$
corresponds the $\beta$-period vector $U^{\pm}_l$:
$$
(U^{\pm}_l)^j = \frac{1}{2\pi i } \int_{\beta_j} \eta_l^{\pm}.
$$

There exist effective divisors $Q_1 + \dots + Q_g$ and $R_1 + \dots + R_g$
such that there are linear equivalences
$$
D = P_1 + \dots + P_{g+1} \ \ \sim  \ \
Q_1 + \dots + Q_g + \infty_-,
$$ 
$$
D = P_1 + \dots + P_{g+1} \ \ \sim  \ \
R_1 + \dots + R_g + \infty_+.
$$
Put $Q_{g+1} =\infty_-$ and $R_{g+1} = \infty_+$
and denote $(A(Q_1) + \dots + A(Q_g))$ and 
$(A(R_1) + \dots A(R_g))$ by $A(Q)$ and $A(R)$, respectively. 

Denote by $\delta$ the Riemann constants vector defined as
follows: for generic $u \in J(\Gamma)$ the multi-valued function
$\vartheta(A(P) - u)$ vanishes exactly at $g$ points $S_1, \dots, S_g$
such that $u + \delta = A(S_1) + \dots + A(S_g)$. 
We also take $\varepsilon$, an odd half-period of 
$\vartheta$ , i.e., $\vartheta(\varepsilon) = 0$ and
$2\varepsilon \equiv 0$ on $J(\Gamma)$.

\bp
The function $\psi$ from Proposition 1 takes the form
$$   
\psi_1(z,\bar{z},P) = 
\exp{\left(z(\int_{P_0}^P \eta^+_1 - a^+_1) + 
\bar{z}(\int_{P_0}^P \eta^-_1 - b^-_1)
\right)}
\cdot 
$$
\be
\frac{\vartheta(A(P) + F_1(z,\bar{z}))}
{\vartheta(A(P) + \delta - A(Q))} \cdot 
\frac{\vartheta(A(\infty_+) + \delta - A(Q))}
{\vartheta(A(\infty_+) + F_1(z,\bar{z}))} \cdot
\label{29}
\ee
$$
\cdot
\frac{\prod^{g+1} \vartheta(\varepsilon + A(P) - A(Q_j))
\cdot \vartheta(\varepsilon + A(\infty_+) - A(P_j))}
{\prod^{g+1} \vartheta(\varepsilon + A(P) - A(P_j))
\cdot \vartheta(\varepsilon + A(\infty_+) - A(Q_j))},
$$
and
$$   
\psi_2(z,\bar{z},P) = 
\exp{\left(z\int_{P_0}^P 
(\eta^+_1 - b^+_1) + \bar{z}(\int_{P_0}^P \eta^-_1 - a^-_1)
\right)}
\cdot
$$
\be 
\frac{\vartheta(A(P) + F_2(z,\bar{z}))}
{\vartheta(A(P) + \delta - A(R))} \cdot 
\frac{\vartheta(A(\infty_-) + \delta - A(R))}
{\vartheta(A(\infty_-) + F_2(z,\bar{z}))} \cdot
\label{30}
\ee
$$
\cdot
\frac{\prod^{g+1} \vartheta(\varepsilon + A(P) - A(R_j))
\cdot \vartheta(\varepsilon + A(\infty_-) - A(P_j))}
{\prod^{g+1} \vartheta(\varepsilon + A(P) - A(P_j))
\cdot \vartheta(\varepsilon + A(\infty_-) - A(R_j))},
$$ 
where the constants $a^{\pm}_1$ and $b^{\pm}_1$ are defined as follows 
\be
\int_{P_0}^P \eta^{\pm}_1 - a^{\pm}_1 = k_{\pm} + O(k^{-1}_{\pm})
\ \mbox{near $\infty_{\pm}$}
\ \ \mbox{and} \ \ 
\int_{P_0}^P \eta^{\pm}_1 - b^{\pm}_1 = O(k^{-1}_{\mp})
\ \ \mbox{near $\infty_{\mp}$}
\label{31}
\ee
and
$$
F_1(z,\bar{z}) = U^+_1 z + U^-_1 \bar{z} + \delta - A(Q),
$$
$$
F_2(z,\bar{z}) = U^+_1 z + U^-_1 \bar{z} + \delta - A(R).
$$
The potentials $U$ and $V$ take the form
\be
U = c_1 \exp{(z(a^+_1-b^+_1)+\bar{z}(b^-_1-a^-_1))}
\frac{\vartheta(A(\infty_+)+ F_2(z,\bar{z}))}
{\vartheta(A(\infty_-)+ F_2(z,\bar{z}))}
\label{32}
\ee
with
$$
c_1 = -  
\frac{\prod^{g+1}\vartheta(\varepsilon + A(\infty_-) - A(P_j))}
{\prod^{g+1} \vartheta(\varepsilon + A(\infty_+) - A(P_j))
\cdot \vartheta(\varepsilon + A(\infty_-) - A(R_j))}
\cdot
$$
$$
\prod^{g}_{j=1} \vartheta(\varepsilon + A(\infty_+) - A(R_j)) 
\cdot
\sum (U^+_1)^j \frac{\partial \vartheta (\varepsilon)}{\partial u^j}
$$
and
\be
V = c_2 \exp{(z(b^+_1-a^+_1)+\bar{z}(a^-_1-b^-_1))}
\frac{\vartheta(A(\infty_-)+ F_1(z,\bar{z}))}
{\vartheta(A(\infty_+)+ F_1(z,\bar{z}))}
\label{33}
\ee
with 
$$
c_2 =   
\frac{\prod^{g+1}\vartheta(\varepsilon + A(\infty_+) - A(P_j))}
{\prod^{g+1} \vartheta(\varepsilon + A(\infty_-) - A(P_j))
\cdot \vartheta(\varepsilon + A(\infty_+) - A(Q_j))}
\cdot
$$
$$
\prod^{g}_{j=1} \vartheta(\varepsilon + A(\infty_-) - A(Q_j)) 
\cdot
\sum (U^-_1)^j \frac{\partial \vartheta (\varepsilon)}{\partial u^j}.
$$
\ep

Here we assume 
that the paths joining $P_0$ and small neighborhoods of infinities 
and the paths coming in definition of $A(P)$ and $A(\infty_{\pm})$ 
are the same and for some homotopy classes of them the expansions
(\ref{31}) hold.

A proof of Proposition 2.

The formulas for $\psi$ are verified by using the periodicity
properties of theta functions and the Riemann theorem on zeroes of
theta functions (\cite{Fay}).

Derive (\ref{32}) and (\ref{33}).
For instance, near $\infty_-$ the function $\psi_1$ decomposes
into the product $\vartheta(\varepsilon + A(P) - A(\infty_-)) 
\cdot H(z,\bar{z},P) \cdot \exp{(k_-\bar{z})}$ and it is known that
$\partial  A(P)/\partial k^{-1}_- = U^-_1$ at $\infty_-$. Hence, we have
$$
\xi^-_{11} = H(z,\bar{z},\infty_-) \cdot 
\sum (U^-_1)^j \frac{\partial\vartheta(\varepsilon)}{\partial u^j}.
$$

This proves Proposition 2.

\bp
Let the spectral data $(\Gamma, \infty_{\pm}, k_{\pm}, D)$ of a finite-zone 
operator $L$ (\ref{25}) (see Proposition 1) satisfy the following conditions

i) there exists a holomorphic involution $\sigma: \Gamma \rightarrow \Gamma$
such that

i1) $\sigma(\infty_{\pm}) = \infty_{\pm}$ and $\sigma(k_{\pm}) = - k_{\pm}$;

i2) there exists a meromorphic differential $\omega$ on $\Gamma$ with 
zeroes at $D + \sigma(D)$ and with two poles at $\infty_{\pm}$
with the principal parts $(\pm k^2_{\pm} + O(k^{-1}_{\pm}))d k^{-1}_{\pm}$;

ii) there exists an antiholomorphic involution 
$\tau:\Gamma \rightarrow \Gamma$ such that 

ii1) $\tau(\infty_{\pm}) = \infty_{\mp}, \tau(k_{\pm}) = -\bar{k}_{\mp}$;

ii2) there exists a meromorphic differential $\tilde{\omega}$ on $\Gamma$ with 
zeroes at $D + \tau(D)$ and with two poles at $\infty_{\pm}$
with the principal parts $(k^2_{\pm} + O(k^{-1}_{\pm}))d k^{-1}_{\pm}$.

Then $L$ takes the form (\ref{2}) with a real potential $U$, i.e., 
$U = V = \bar{U}$.
\ep 

A proof of Proposition 3.

Consider the meromorphic
differential $\psi_1(P)\psi_2(\sigma(P))\omega$
with poles only at $\infty_{\pm}$. The sum of its residue equals
$\xi^+_{21} + \xi^-_{11}$. Since it vanishes,
by (\ref{28}), this implies $U = V$.

Now consider the differentials
$\psi_1(P)\overline{\psi_1(\tau(P))}\tilde{\omega}$
and $\psi_2(P)\overline{\psi_2(\tau(P))}\tilde{\omega}$
and, computing the sums of their residues as above, we conclude
$U = \bar{U}$ and $V=\bar{V}$.

This proves the proposition.

\bp
1) Given a spectral data $(\Gamma, \infty_{\pm}, k_{\pm}, D)$ meeting the 
conditions of Proposition 3, for any $n$-tuple of points
$Q_1, \dots, Q_n \in \Gamma \setminus \{\infty_{\pm}, P_1, \dots, P_g\}$
and any $n$-tuple of constants $a_1, \dots, a_n \in \C$
the vector function $\psi(z,\bar{z}) = a_1 \psi(z,\bar{z},Q_1) + \dots +
a_n \psi(z,\bar{z},Q_n)$ defines via (\ref{3}) a surface immersed into
$\R^3$.

2) If $\psi(z+\gamma) = \pm \psi(z)$ for $z \in \C$ and 
$\gamma \in \Lambda$, with $\Lambda \subset \C$ a rank two lattice,
then the immersion converts into an immersion of the torus 
$\C/\Lambda$ if and only if
\be
\int_{\C/\Lambda} \psi^2_1 dz \wedge d\bar{z} =
\int_{\C/\Lambda} \psi^2_2 dz \wedge d\bar{z} =
\int_{\C/\Lambda} \psi_1 \bar{\psi}_2 dz \wedge d\bar{z} = 0.
\label{34}
\ee
\ep

A proof of Proposition 4.

The first statement is evident. The second statement
is quite clear for tori of revolution 
(see \cite{T2}) and has been extended onto the general case
by M. Schmidt. 

Let $\hat{\gamma}$ and $\tilde{\gamma}$ be generators of 
$\Lambda$ with the basis $(\hat{\gamma},\tilde{\gamma})$ is positively
oriented. By (\ref{3}), $X^1 + i X^2 = X^1 + i \bar{X}^2 = 
i\int (\bar{\psi}^2_2 dz - \bar{\psi}^2_1 d \bar{z})$.
We have $d( x(\bar{\psi}^2_2 dz - \bar{\psi}^2_1 d \bar{z})) = 
i(\bar{\psi}^2_1 + \bar{\psi}^2_2) dx \wedge dy$ and
$d(y(\bar{\psi}^2_2 dz - \bar{\psi}^2_1 d \bar{z})) =
(\bar{\psi}^2_1 - \bar{\psi}^2_2) dx \wedge dy$.
The Stokes theorem implies that
$$
\int_{\C/\Lambda} 
(\bar{\psi}^2_1 + \bar{\psi}^2_2) dz' \wedge d\bar{z}'=
-2(\mbox{Re}\, \hat{\gamma} \int_{z}^{z+\tilde{\gamma}}
(\bar{\psi}^2_2 dz' - \bar{\psi}^2_1 d \bar{z}') 
-
\mbox{Re}\,\tilde{\gamma} \int_{z}^{z+\hat{\gamma}}
(\bar{\psi}^2_2 dz' - \bar{\psi}^2_1 d \bar{z}')
)
$$
and
$$
\int_{\C/\Lambda} 
(\bar{\psi}^2_1 - \bar{\psi}^2_2) dz' \wedge d\bar{z}' =
2i(\mbox{Im}\,\hat{\gamma} \int_{z}^{z+\tilde{\gamma}}
(\bar{\psi}^2_2 dz' - \bar{\psi}^2_1 d \bar{z}')
-
\mbox{Im}\,\tilde{\gamma} \int_{z}^{z+\hat{\gamma}}
(\bar{\psi}^2_2 dz' - \bar{\psi}^2_1 d \bar{z}')).
$$
Therefore $X^1 + iX^2$ is $\Lambda$-periodic if and only if 
the first two equalities from (\ref{34}) hold. The equivalence of the 
last one to $\Lambda$-periodicity of $X^3$ is proven in the same manner. 
This proves the proposition
\footnote{The constructed surfaces may have singularities 
which are exactly at points where
$|\psi_1|^2 + |\psi_2|^2 =0$. 
For $C^2$-regular surfaces the periodicity conditions (\ref{34}) are 
quite perspequitive because, by (\ref{7}) and (\ref{8}),  
$\psi_1^2 = -\partial(X^2+iX^1)$,
$\psi_2^2 = \bar{\partial}(X^2+iX^1)$, and $\psi_1\bar{\psi}_2 = 
\partial X^3$ globally.}.

\bc 
{\bf \S 5. Finite-zone solutions to the modified Novikov--Veselov equations}
\ec

The modified Novikov--Veselov (mNV) equations are related to 
the operator ${\cal D}$ (\ref{3}) and take the form
\be
{\cal D}_t = {\cal D} A + B {\cal D} 
\label{35}
\ee
(the Manakov triple). 
The deformations of $U$ generate the deformations of $\psi$, 
the zero-eigenfunction of ${\cal D}$, of the form
\be
\psi_t = A\psi.
\label{36}
\ee
These equations 
had been introduced by Bogdanov in \cite{Bg} and it is
an observation of Konopelchenko 
that if a surface is 
defined  by (\ref{3}) then (\ref{35}) generates via (\ref{36}) a local
deformation of the surface (\cite{Kon}).  

The first equation of this hierarchy is
$$
U_t = (U_{zzz} + 3 U_zV + \frac{3}{2}UV_z) + 
(U_{\bar{z}\bar{z}\bar{z}} + 3U_{\bar{z}}\bar{V} + 
\frac{3}{2}U\bar{V}_{\bar{z}})
$$
with $V_{\bar{z}} = (U^2)_{z}$.
If $U$ depends only on one spatial variable then these equations
reduce to the modified
Korteweg--de Vries equations. 

\bp
Let  the spectral data $(\Gamma, \infty_{\pm}, k_{\pm}, D)$ are as in 
Proposition 3. Define the constants $a^{\pm}_{l}$ and $b^{\pm}_l$ 
by the following analogs of (\ref{31})  
$$ 
\int_{P_0}^P \eta^{\pm}_l - a^{\pm}_l = k^l_{\pm} + O(k^{-1}_{\pm})
\ \mbox{near $\infty_{\pm}$} \ \ \mbox{and} \ \ 
\int_{P_0}^P \eta^{\pm}_l - b^{\pm}_l = O(k^{-1}_{\mp})
\ \mbox{near $\infty_{\mp}$}
$$
and define $F_{1l}(z,\bar{z},t_l)$ and $F_{2l}(z,\bar{z},t_l)$
by
$$
F_{1l}(z,\bar{z},t_l) = U^+_1z+U^-_1\bar{z}+ (U^+_{2l+1}+U^-_{2l+1})t_l +
\delta - A(Q),
$$
$$
F_{2l}(z,\bar{z},t_l) = U^+_1z+U^-_1\bar{z}+ (U^+_{2l+1}+U^-_{2l+1})t_l +
\delta - A(R).
$$
Let $\hat{\psi}$ be a vector function obtained from (\ref{29}) and (\ref{30})
by replacing $F_1$ and $F_2$ by $F_{1l}$ and $F_{2l}$
and by adding 
$t_l (\int_{P_)}^P (\eta^+_{2l+1} + \eta^-_{2l+1}) - (a^+_{2l+1}+
b^-_{2l+1}))$ 
and
$t_l (\int_{P_)}^P (\eta^+_{2l+1} + \eta^-_{2l+1}) - (a^-_{2l+1}+
b^+_{2l+1}))$ 
to the arguments of exponents in (\ref{29}) and (\ref{30}).
Then ${\cal D}\hat{\psi} = 0$ where $U(z,\bar{z},t_l)$
is constructed by $\hat{\psi}$ via (\ref{28}) and satisfies
the $l$-th equation of the mNV hierarchy. 
\ep
 
The formula for $U(z,\bar{z},t_l)$ is derived by  
the same substitutions. 

The proof of this proposition is 
analogous to one for the Novikov--Veselov equation (\cite{VN2})
as well as Proposition 3 is an analog of the Novikov--Veselov theorem 
distinguishing potential finite-zone Schr\"odinger operators (\cite{VN1}).
But there is one important difference. 
The form $\omega$ from Proposition 3 is not invariant under $\sigma$ but
antiinvariant, i.e., $\sigma^*(\omega) = -\omega$, and this implies 
the following. 
Finite-zone potentials  can be 
written in terms of Prym theta functions 
of the covering $\Gamma \rightarrow \Gamma/\sigma$ 
but we cannot control the topological type of
involution $\sigma$ and 
cannot conclude that for nonsingular $\Gamma$ it has
only two fixed points. We can only say that
\be
\mbox{genus} (\Gamma) - \mbox{genus}(\Gamma/\sigma) = 
\dim \mbox{Prym}(\Gamma,\sigma) \geq [\frac{\mbox{genus}(\Gamma)}{2}].
\label{37}
\ee

\bc
{\bf \S 6. The Willmore functional}
\ec

The following conjecture of Willmore is well-known 

{\sl For tori immersed into $\R^3$ the minimum of the Willmore functional
$$
{\cal W}(\Sigma) = \int_{\Sigma} H^2 d\mu,
$$
with $d\mu$ an induced Liouville measure, equals $2\pi^2$ and is attained
on the Clifford torus and its images under conformal transformations of
$\R^3$.}

The Clifford torus is obtained by revolving of a circle of radius 
$1$ around the axis lying in the same plane as the circle at distance
$\sqrt{2}$ from the circle center.

Until now this conjecture is proved only in some particular cases
(see the survey of them in \cite{T1}). 

By (\ref{4}), for a torus $\Sigma \subset \R^3$ the simplest
Kruskal integral 
$$
4 \int_{\C/\Lambda} U^2(z,\bar{z}) dx dy, 
$$
where $U : \C/\Lambda \rightarrow \C$ is the potential of its Weierstrass
representation, coincides with ${\cal W}(\Sigma)$ (\cite{T1}).
This leads us to the conjecture that 

{\sl for fixed conformal classes of tori
the minima of ${\cal W}$ are attained on tori with minimal
spectral genus.}

This fits into the Willmore conjecture because the spectral genus of
the Clifford torus equals zero (this is derived from its Weierstrass
representation found in \cite{T2}). Moreover, it is difficult to 
imagine that the minimum of such variational problem is degenerated. 
In \cite{T1} it was conjectured
that the minima of this functional are stationary with respect to 
the deformations generated by the first mNV equation
for which it is shown in (\cite{T1}) that this deformation preserves tori.
We may extend that as follows

{\sl for all equations of the mNV  hierarchy 
the minima of ${\cal W}$ for fixed conformal classes
are stationary with respect to induced deformations.}

This  implies that the Prym variety of $\Gamma \rightarrow \Gamma/\sigma$
would be one-dimensional and the mNV deformations
reduce to translations of tori along themselves.    
By (\ref{37}), this implies $\mbox{genus}(\Gamma) \leq 3$.
For surfaces of genus $3$ the dimension of the Prym variety equals
$3$ if $\sigma$ is a hyperelliptic involution or $2$
if $\sigma$ has $4$ fixed points (the case when this dimension equals $1$
correspond to an involution without fixed points).
Hence, the last conjecture implies that for minima $\mbox{genus}(\Gamma)
\leq 2$.

The Weierstrass representation gives a physical 
explanation for lower bounds for ${\cal W}$: 
it is clear that for small perturbations of the zero-potential
$U=0$ the surfaces constructed in terms of zero-eigenfunctions of
(\ref{2}) do not convert into tori and, since for $U$ the Willmore functional
is its $L_2$-norm, the lower bound reflects how big a perturbation has to be
to get planes converted into tori.

{\bf Final remark.}

This work is supported by the Russian Foundation 
for Basic Researches
(grant 96-15-96877) and by SFB 288. The part of this work has been done
during the author's stay at IHES (January 1997).
The author thanks P. Grinevich, U. Pinkall, and M. Schmidt for helpful
conversations.

\end{document}